\documentclass[11pt]{article}
\usepackage{latexsym}
\usepackage{amsfonts}
\usepackage{amssymb}
\usepackage{mathrsfs}

\usepackage{epic}
\usepackage{eepic}
\usepackage{hyperref}

\newcommand{\be}{\begin{equation}}
\newcommand{\ee}{\end{equation}}
\newcommand{\bea}{\begin{eqnarray}}
\newcommand{\eea}{\end{eqnarray}}
\newcommand{\bean}{\begin{eqnarray*}}
\newcommand{\eean}{\end{eqnarray*}}
\newcommand{\brray}{\begin{array}}
\newcommand{\erray}{\end{array}}

\newcommand{\newsection}[1]{\setcounter{equation}{0}
\setcounter{dfn}{0}
\section{#1}}

\newtheorem{dfn}{Definition}[section]
\newtheorem{thm}[dfn]{Theorem}
\newtheorem{lmma}[dfn]{Lemma}
\newtheorem{ppsn}[dfn]{Proposition}
\newtheorem{crlre}[dfn]{Corollary}
\newtheorem{xmpl}[dfn]{Example}
\newtheorem{rmrk}[dfn]{Remark}

\newcommand{\bdfn}{\begin{dfn}}
\newcommand{\bthm}{\begin{thm}}
\newcommand{\blmma}{\begin{lmma}}
\newcommand{\bppsn}{\begin{ppsn}}
\newcommand{\bcrlre}{\begin{crlre}}
\newcommand{\bxmpl}{\begin{xmpl}}
\newcommand{\brmrk}{\begin{rmrk}}

\newcommand{\edfn}{\end{dfn}}
\newcommand{\ethm}{\end{thm}}
\newcommand{\elmma}{\end{lmma}}
\newcommand{\eppsn}{\end{ppsn}}
\newcommand{\ecrlre}{\end{crlre}}
\newcommand{\exmpl}{\end{xmpl}}
\newcommand{\ermrk}{\end{rmrk}}

\newcommand{\bbc}{\mathbb{C}}
\newcommand{\bbz}{\mathbb{Z}}
\newcommand{\bbn}{\mathbb{N}}

\newcommand{\scrf}{\mathscr{F}}

\newcommand{\cla}{\mathcal{A}}
\newcommand{\clb}{\mathcal{B}}
\newcommand{\clh}{\mathcal{H}}

\newcommand{\cll}{\mathcal{L}}

\newcommand{\halpha}{\widehat{\alpha}}
\newcommand{\hbeta}{\widehat{\beta}}
\newcommand{\hgamma}{\widehat{\gamma}}

\newcommand{\prf}{\noindent{\it Proof\/}: }

\def \qed { \mbox{}\hfill
$\Box$\vspace{1ex}}

\newcommand{\half}{\frac{1}{2}}

\begin{document}


\author{{\sc Partha Sarathi Chakraborty} and
{\sc Arupkumar Pal}}
\title{Equivariant spectral triples and Poincar\'e duality for $SU_q(2)$}
\maketitle
   \begin{abstract}
Let $\cla$ be the $C^*$-algebra associated with $SU_q(2)$, $\pi$ be
the representation by left multiplication on the $L_2$ space of the Haar
state and let $D$ be the equivariant Dirac operator for this representation
constructed by the authors earlier.  We prove
in this article that there is no operator other than the scalars
in the commutant $\pi(\cla)'$ that has bounded commutator with $D$.
This  implies that the equivariant spectral triple under consideration
 does not admit a rational Poincar\'e dual in the sense of Moscovici, which in particular
means that this spectral triple
 does not extend to a $K$-homology fundamental class for $SU_q(2)$.
We also show that a minor modification of this equivariant spectral triple
gives a fundamental class and thus implements Poincar\'e duality. 
 \end{abstract}
{\bf AMS Subject Classification No.:} {\large 58}B{\large 34}, {\large 46}L{\large 87}, {\large
  19}K{\large 33}

\newsection{Introduction}
 In noncommutative geometry (NCG), spaces are described by a
triple $(\cla,\clh,D)$, where $\mathcal{A}$ is a $*$-algebra closed under
holomorphic functional calculus acting on a complex separable
Hilbert space $\mathcal{H}$ and $D$ is an unbounded self-adjoint
operator with compact resolvent that has bounded commutators
with elements from the algebra $\mathcal{A}$. Such a triple is called a spectral triple.
In this spectral point of view,  one requires $D$ to
be nontrivial in the sense that the
associated Kasparov module should give a nontrivial element in
$K$-homology.
One can also formulate the notion of Poincar\'e duality
in this context.
A pair of separable $C^∗$-algebras $(A,B)$ is said
to be a \textbf{Poincar\'e dual pair}  if there exist a class 
$\Delta \in KK(A\otimes B,\bbc)$ and a class 
$\delta\in KK(\bbc,A\otimes B)$ with the properties
$\delta \otimes_B \Delta = id_A\in KK(A,A)$ and 
$\delta \otimes_A \Delta = id_B\in KK(B,B)$.
The element $\Delta$ is called a $K$-homology fundamental class for 
the pair $(A,B)$.
Poincar\'e duality is said to hold for a separable $C^*$-algebra $A$
if there is a $K$-homology fundamental class for the pair $(A,A)$.
See section~4, chapter~6 in \cite{co1} for a detailed formulation, and~\cite{m}
for an interesting application. 

The existence of a fundamental class can often be 
deduced from abstract $KK$-theory arguments,
using the properties of the $C^*$-algebra in question.
But more interesting from the point of view of 
noncommutative geometry is an 
explicit geometric realization of this fundamental
class, with possibly other nice features. 
An explicit geometric 
realization of a $K$-homology class is given by
a spectral triple. Suppose we have a spectral triple
$(\clh, \pi, D)$ for a $C^*$-algebra $A$, where $\pi$ is faithful. 
If there is another faithful representation $\pi'$ of $A$ on $\clh$
such that 
\begin{enumerate}
\item 
$\pi'$ and $\pi$ commute, 
\item 
$(\clh,\pi',D)$ is a spectral triple for $A$,
and 
\item 
$(\clh, \pi\otimes\pi', D)$ gives a $K$-homology fundamental
class for $A\otimes A$, 
\end{enumerate}
then we say that the spectral triple
$(\clh, \pi,D)$ extends to a $K$-homology fundamental class.
If one replaces condition~3 above with a slightly weaker condition,
then one says that the spectral triples $(\clh, \pi, D)$ and 
$(\clh,\pi',D)$ are rational Poincar\'e duals (see~\cite{m} for
this notion).
In an earlier paper (\cite{c-p}), the authors constructed
an equivariant spectral triple for the quantum $SU(2)$ group
that was later analysed further by Connes in~\cite{co3}.
It is natural to ask whether the triple gives rise to a fundamental
class for $SU_q(2)$. This is what we try to answer in this paper.

Let $h$ be the Haar state for the quantum $SU(2)$ group
and let $\pi$ be the representation of $C(SU_q(2))$ on $L_2(h)$
by left multiplication. 
In section~2, we make a detailed analysis of
the operators $\alpha$ and $\beta$ on $L_2(h)$. We also introduce
and study two operators  
$\halpha$ and $\hbeta$ that are compact perturbations
of $\alpha$ and $\beta$ respectively and obey the same
commutation relations as $\alpha$ and $\beta$.
These play an important role in the proof of the main
result in section~4.
In section~3, we compute the modular conjugation operator associated
with the Haar state.  This helps us describe elements
of the commutator in terms of elements of the strong closure
of $\pi(C(SU_q(2)))$.
Denote by $D$ the equivariant Dirac operator constructed
by the authors in~\cite{c-p}.
In section~4, we prove that there is no operator other than the scalars
in the commutant of $\pi(C(SU_q(2)))$, that has bounded commutator with $D$.
An important consequence of this is that the equivariant spectral
triple does not give a $K$-homology fundamental class for $SU_q(2)$.
In the final section, we show that Poincar\`e duality holds
for $SU_q(2)$. We also give an explicit construction of a spectral
triple that gives a fundamental class for $SU_q(2)$.

\newsection{Closer look at the $L_2$ space}
In what follows, we will be concerned with  
the quantum $SU(2)$ group, the spectral triple under
consideration being the equivariant spectral triple
constructed by the authors in \cite{c-p}.
To fix notation, let us recall a few things
from that paper.
Let $q$ be a real number in the interval $(0,1)$.
Let $\cla$ denote the $C^*$-algebra of
continuous functions on $SU_q(2)$,
which is the universal $C^*$-algebra
generated by two elements $\alpha$ and $\beta$
subject to the relations
\be\label{suq2comm}
\alpha^*\alpha+\beta^*\beta=I=\alpha\alpha^*+q^2\beta\beta^*,
\quad
\alpha\beta-q\beta\alpha=0=\alpha\beta^*-q\beta^*\alpha,
\quad
\beta^*\beta=\beta\beta^*
\ee
Let $h$ denote the Haar state on $\cla$ and
let $\pi :\cla\rightarrow\cll(L_2(h))$ be
the representation given by left multiplication
by elements in $\cla$.
We will often identify an element $a\in\cla$ with $\pi(a)$.
$\alpha_r$ and $\beta_r$ will stand for
$\alpha^r$ and $\beta^r$ respectively if $r\geq 0$, and
for $(\alpha^*)^{-r}$ and $(\beta^*)^{-r}$ if $r<0$.
Let $D$ be the operator given by
$D:e^{(n)}_{ij}\mapsto d(n,i)e^{(n)}_{ij}$, where
\be \label{genericd}
d(n,i)=\cases{2n+1 & if $n\neq i$,\cr
                      -(2n+1) & if $n=i$.}
\ee
Then $(L_2(h),\pi, D)$ is an odd equivariant spectral triple
of dimension 3 and with nontrivial $K$-homology class.

Our objective is to study commutators of the form $[D,T']$
with $T'$ coming from the commutant $(\pi(\cla))'$.
Any such $T'$ can be written as $JTJ$ where $J$ is the modular conjugation
operator associated with the Haar state and $T$ comes from
the strong closure $(\pi(\cla))''$ of $\pi(\cla)$.
With this in mind, in this section we study the structures of the
operators that constitute $(\pi(\cla))''$.
Recall (cf.\ \cite{c-p}) that
$L_2(h)$ has a natural orthonormal basis
$\{e^{(n)}_{ij}: n\in\half\bbn, i,j=-n,-n+1,\ldots,n\}$,
and the left multiplication operators in this basis are given by
\bea
\alpha: e^{(n)}_{ij}  &\mapsto& a_+(n,i,j) e^{(n+\half)}_{i-\half ,j-\half }
       + a_-(n,i,j)  e^{(n-\half )}_{i-\half ,j-\half },\label{alpha}\\
\beta:e^{(n)}_{ij}  &\mapsto& b_+(n,i,j)  e^{(n+\half )}_{i+\half ,j-\half }
       + b_-(n,i,j)  e^{(n-\half )}_{i+\half ,j-\half },\label{beta}
\eea
where
\bean
a_+(n,i,j)  & =& \Bigl(q^{2(n+i)+2(n+j)+2}
 \frac{(1-q^{2n-2j+2})(1-q^{2n-2i+2})}{(1-q^{4n+2})(1-q^{4n+4})}
                          \Bigr)^\half,\label{aplus}\\
a_-(n,i,j)&=&\Bigl(\frac{(1-q^{2n+2j})(1-q^{2n+2i})}
                                   {(1-q^{4n})(1-q^{4n+2})}\Bigr)^\half,\label{aminus}\\
b_+(n,i,j)&=& - \Bigl(q^{2(n+j)}\frac{(1-q^{2n-2j+2})(1-q^{2n+2i+2})}
                              {(1-q^{4n+2})(1-q^{4n+4})}\Bigr)^\half,\label{bplus}\\
b_-(n,i,j)&=&\Bigl(q^{2(n+i)}\frac{(1-q^{2n+2j})(1-q^{2n-2i})}{(1-q^{4n})
   (1-q^{4n+2})}\Bigr)^\half.\label{bminus}
\eean

We will also need the following operators on $L_2(h)$:
\bea
\halpha: e^{(n)}_{ij}  &\mapsto&
 \hat{a}_+(n,i,j) e^{(n+\half)}_{i-\half ,j-\half }
  + \hat{a}_-(n,i,j)  e^{(n-\half )}_{i-\half ,j-\half },\label{halpha}\\
\hbeta:e^{(n)}_{ij}  &\mapsto& \hat{b}_+(n,i,j)  e^{(n+\half )}_{i+\half ,j-\half }
      + \hat{b}_-(n,i,j)  e^{(n-\half )}_{i+\half ,j-\half },\label{hbeta}
\eea
where
\bean
\hat{a}_+(n,i,j)  & =& q^{2n+i+j+1},\label{caplus}\\
\hat{a}_-(n,i,j)&=&(1-q^{2n+2i})^\half(1-q^{2n+2j})^\half,\label{caminus}\\
\hat{b}_+(n,i,j)&=& - q^{n+j}(1-q^{2n+2i+2})^\half,\label{cbplus}\\
\hat{b}_-(n,i,j)&=&q^{n+i}(1-q^{2n+2j})^\half.\label{cbminus}
\eean
It is easy to see that $\halpha$ and $\hbeta$ are compact perturbations of
$\alpha$ and $\beta$ respectively.

We will now decompose the space $L_2(h)$ as a direct sum
of smaller subspaces, and study the behaviour of the above operators
with respect to this decomposition.
Note that the set $\Lambda=\{(n,i,j): n\in \half\bbn, i,j=-n,-n+1,\ldots,n\}$
parametrizes the canonical orthonormal basis for $L_2(h)$.
For each $n\in\half\bbz$, denote by $\Lambda_n$ the minimal subset of
$\Lambda$ containing the point $(|n|,-n,-n)$ and closed under the
translations
\[
(a,b,c)\mapsto (a+\half,b+\half, c-\half),\quad
(a,b,c)\mapsto (a+\half,b-\half, c+\half).
\]
For $n,k\in\half\bbz$, denote by $\Lambda_{nk}$ the minimal subset of $\Lambda_n$ that contains
$(|n|+|k|, -n+k, -n-k)$ and is closed under the translation
$(a,b,c)\mapsto (a+1,b,c)$.
Thus all the $\Lambda_{nk}$'s are disjoint,
$\Lambda=\cup_n \Lambda_n$,
$\Lambda_n=\cup_{k}\Lambda_{nk}$.
The following diagram will make it easier to visualize these sets.
Represent the lattice $\Lambda$ as a pyramid, where the vertical axis
is the $n$-axis, the top vertex is on the plane $n=0$ and $n$ increases
downwards.
Then $\Lambda_n$ are precisely the vertical cross-sections parallel to
the plane $ABD$. 
The ones that intersect the triangle $BCD$ correspond to
nonnegative values of $n$ and the ones that intersect the triangle $BDE$
correspond to nonpositive values of $n$.
In particular, $\Lambda_0$ is the cross-section given by the plane $ABD$.
Similarly
$\Lambda_{nk}$ are vertical lines in the plane $\Lambda_n$.
The lines that intersect the triangle $CDE$ correspond to nonnegative
values of $k$ and the lines that intersect the triangle $BCE$
correspond to nonpositive values of $k$.
In particular, $\Lambda_{n0}$ are the lines that intersect the line $CE$.

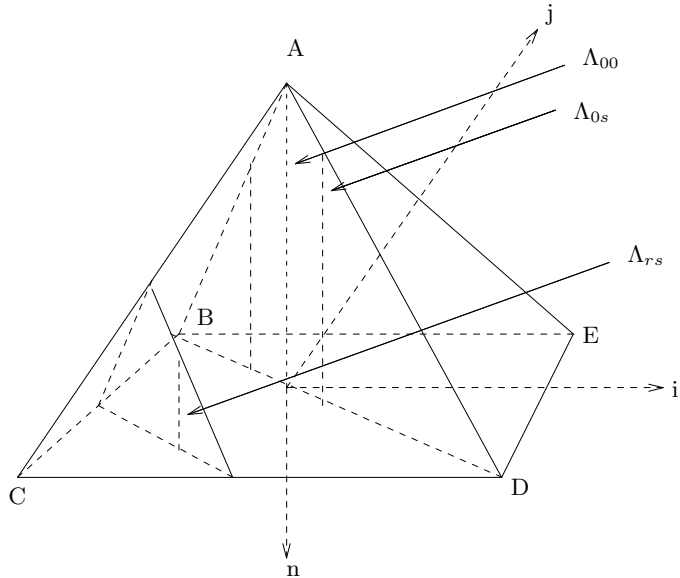
\begin{center}
\setlength{\unitlength}{0.00062500in}
\begingroup\makeatletter\ifx\SetFigFont\undefined%
\gdef\SetFigFont#1#2#3#4#5{%
  \reset@font\fontsize{#1}{#2pt}%
  \fontfamily{#3}\fontseries{#4}\fontshape{#5}%
  \selectfont}%
\fi\endgroup%
{\renewcommand{\dashlinestretch}{30}
\begin{picture}(5758,5200)(0,1000)
\drawline(2487,5412)(237,2112)
\drawline(2487,5412)(4287,2112)
\drawline(237,2112)(4287,2112)
\dashline{60.000}(1587,3312)(2487,5412)
\drawline(2487,5412)(4887,3312)(4287,2112)
\dashline{60.000}(1587,3312)(4887,3312)
\drawline(1362,3687)(2037,2112)
\drawline(12,12)(12,12)
\drawline(12,12)(12,12)
\drawline(2964.745,4580.873)(2862.000,4512.000)(2985.068,4524.420)
\texture{44555555 55aaaaaa aa555555 55aaaaaa aa555555 55aaaaaa aa555555 55aaaaaa
        aa555555 55aaaaaa aa555555 55aaaaaa aa555555 55aaaaaa aa555555 55aaaaaa
        aa555555 55aaaaaa aa555555 55aaaaaa aa555555 55aaaaaa aa555555 55aaaaaa
        aa555555 55aaaaaa aa555555 55aaaaaa aa555555 55aaaaaa aa555555 55aaaaaa }
\drawline(2862,4512)(4737,5187)
\drawline(2862,4512)(4737,5187)
\drawline(1764.641,2706.028)(1662.000,2637.000)(1785.049,2649.605)
\drawline(1662,2637)(5187,3912)
\drawline(1662,2637)(5187,3912)
\dashline{60.000}(1587,3087)(1587,2337)
\dashline{60.000}(2037,2112)(912,2712)(1362,3762)
\dashline{60.000}(2787,4812)(2787,2712)
\dashline{45.000}(2487,5412)(2487,2862)
\dashline{60.000}(2187,4737)(2187,3012)
\dashline{60.000}(4287,2112)(1512,3312)
\dashline{60.000}(237,2112)(1587,3312)
\dashline{60.000}(2487,2862)(5637,2862)
\drawline(5517.000,2832.000)(5637.000,2862.000)(5517.000,2892.000)
\dashline{60.000}(2487,2862)(2487,1437)
\drawline(2457.000,1557.000)(2487.000,1437.000)(2517.000,1557.000)
\dashline{60.000}(2487,2862)(4587,5862)
\drawline(4542.761,5746.488)(4587.000,5862.000)(4493.608,5780.896)
\drawline(4812,5562)(2562,4737)
\drawline(4812,5562)(2562,4737)
\drawline(2664.338,4806.477)(2562.000,4737.000)(2684.993,4750.144)
\put(2487,5637){\makebox(0,0)[lb]{\smash{{{\SetFigFont{9}{10.8}{\rmdefault}{\mddefault}{\updefault}A}}}}}
\put(1737,3387){\makebox(0,0)[lb]{\smash{{{\SetFigFont{9}{10.8}{\rmdefault}{\mddefault}{\updefault}B}}}}}
\put(4362,1962){\makebox(0,0)[lb]{\smash{{{\SetFigFont{9}{10.8}{\rmdefault}{\mddefault}{\updefault}D}}}}}
\put(4887,5112){\makebox(0,0)[lb]{\smash{{{\SetFigFont{9}{10.8}{\rmdefault}{\mddefault}{\updefault}$\Lambda_{0s}$}}}}}
\put(5337,3912){\makebox(0,0)[lb]{\smash{{{\SetFigFont{9}{10.8}{\rmdefault}{\mddefault}{\updefault}$\Lambda_{rs}$}}}}}
\put(2487,1287){\makebox(0,0)[lb]{\smash{{{\SetFigFont{9}{10.8}{\familydefault}{\mddefault}{\updefault}n}}}}}
\put(5712,2787){\makebox(0,0)[lb]{\smash{{{\SetFigFont{9}{10.8}{\familydefault}{\mddefault}{\updefault}i}}}}}
\put(4662,5937){\makebox(0,0)[lb]{\smash{{{\SetFigFont{9}{10.8}{\familydefault}{\mddefault}{\updefault}j}}}}}
\put(4962,5562){\makebox(0,0)[lb]{\smash{{{\SetFigFont{9}{10.8}{\familydefault}{\mddefault}{\updefault}$\Lambda_{00}$}}}}}
\put(162,1887){\makebox(0,0)[lb]{\smash{{{\SetFigFont{9}{10.8}{\familydefault}{\mddefault}{\updefault}C}}}}}
\put(4962,3237){\makebox(0,0)[lb]{\smash{{{\SetFigFont{9}{10.8}{\familydefault}{\mddefault}{\updefault}E}}}}}
\end{picture}
}
\end{center}

Let us also note that the family of maps
$\phi_n:\Lambda_n\rightarrow \Lambda_0$
given by
\be\label{bijections}
\phi_n(a,b,c)=(a-|n|,b+n,c+n)
\ee
give bijections between $\Lambda_n$ and $\Lambda_0$
whose restriction to $\Lambda_{nk}$ yield a bijection
from $\Lambda_{nk}$ to  $\Lambda_{0k}$.

Let $\clh_r$ denote the closed span of $\{e^{(n)}_{ij}: (n,i,j)\in \Lambda_r\}$,
$\clh_{rs}$ denote the closed span of $\{e^{(n)}_{ij}: (n,i,j)\in \Lambda_{rs}\}$,
$P_r$ denote the projection onto $\clh_r$ and $P_{rs}$ denote the projection onto
$\clh_{rs}$.
For an operator $T$, denote by $T_{rs}$ the restriction $P_{rs}TP_{rs}$
of $T$ to $\clh_{rs}$.
Let $U_n$ denote the unitary operator from $\clh_n$ to $\clh_0$
induced by the bijection $\phi_n$.

\bppsn\label{inv1}
Let $A$ stand for $\alpha$ or $\halpha$, and $B$ stand for $\beta$ or $\hbeta$.
Then one has
\bea
P_{n+\half}AP_n=AP_n, && P_{r+\half,s}AP_{rs}=AP_{rs},\label{alphacomm}\\
BP_n=P_nB,&& P_{r,s+\half}BP_{rs}=BP_{rs},\label{betacomm}\\
P_{rs}B^*B=B^*BP_{rs},&& \label{gmcomm}
\eea
where $n,r,s\in\half\bbz$.

Moreover, for all $n\in\half\bbn$, the operators
$U_n\halpha U_n^*$ and $U_n\hbeta U_n^*$ are independent of $n$.
\eppsn
\prf
This is a simple consequence of equations~(\ref{alpha}--\ref{hbeta}).
\qed

\blmma\label{comm}
$\halpha$ and $\hbeta$ satisfy the following commutation relations:
\be
\halpha^*\halpha+\hbeta^*\hbeta=I,\;\halpha\halpha^*
+q^2\hbeta\hbeta^*=I,\;
\halpha\hbeta-q\hbeta\halpha=0,\;\halpha\hbeta^*-q\hbeta^*\halpha=0,\;
\hbeta^*\hbeta=\hbeta\hbeta^*.
\ee
\elmma
\prf
The relations follow by direct computation from the actions
of $\halpha$ and $\hbeta$ given in equations~(\ref{halpha}) and (\ref{hbeta}).
\qed

Following is a simple consequence of the above commutation relations.
\bcrlre \label{temp1}
Let $\gamma=\beta^*\beta$ and $\hgamma=\hbeta^*\hbeta$.
Then $\sigma(\hgamma)=\{q^{2k}:k\in\bbn\}\cup\{0\}=\sigma(\gamma)$,
and $\ker \halpha^*=\{0\}=\ker \alpha^*$.
\ecrlre
Note that the action of $\hgamma$ on the basis vectors are given by
\be\label{gmaction}
\hgamma e^{(n)}_{ij}=
c_+(n,i,j) e^{(n+1)}_{ij}
+c_0(n,i,j)e^{(n)}_{ij}
+c_-(n,i,j) e^{(n-1)}_{ij},
\ee
where
\bean
c_+(n,i,j)  &=& -q^{2n+i+j+1}(1-q^{2n+2i+2})^\half (1-q^{2n+2j+2})^\half,\\
c_0(n,i,j)  &=& (q^{2n+2j}(1-q^{2n+2i})+q^{2n+2i}(1-q^{2n+2j+2}),\\
c_-(n,i,j)  &=& -q^{2n+i+j-1}(1-q^{2n+2i})^\half (1-q^{2n+2j})^\half.
\eean
One can check, using (\ref{gmaction}) and (\ref{gmcomm}),
that $\ker\hgamma=\{0\}$.

\blmma\label{hgamma}
Let $r\in\half\bbn$, and $s\in\half\bbz$.
The restriction $P_{rs}\hgamma P_{rs}$ of $\hgamma$ to $\clh_{rs}$
 is compact and
the spectrum $\sigma(P_{rs}\hgamma P_{rs})$ coincides with $\sigma(\hgamma)$.
\elmma
\prf
Observe that for $r\in\half\bbn$,
$U_r(P_{rs}\hgamma P_{rs})U_r^*=P_{0s}\hgamma P_{0s}$.
So it is enough to prove the statement for $r=0$.

It is easy to see that $P_{0s}\hgamma P_{0s}$ is compact
by using equation~(\ref{gmaction}).
This, along with the second equality in (\ref{betacomm}) and the fact that
$\hbeta$ and $\hgamma$ commute, tells us that
 $\sigma(P_{0s}\hgamma P_{0s})$ is independent of $s$, and
consequently
$\sigma(P_{0s}\hgamma P_{0s})=\sigma(P_0\hgamma P_0)$
and in fact, this is same as the essential spectrum $\sigma_{ess}(P_0\hgamma P_0)$.

Let us next show that $\sigma(P_0\hgamma P_0)=\sigma(\hgamma)$.
Let $K$ be the operator on $\clh_0$, given on the basis
vectors $e^{(n)}_{i,-i}$ as follows:
\be\label{cpt}
K e^{(n)}_{i,-i}=
c_+(n,i,-i) e^{(n+1)}_{i,-i}
+(q^{2n+2|i|}-q^{4n}-q^{4n+2})e^{(n)}_{i,-i}
+c_-(n,i,-i) e^{(n-1)}_{i,-i}.
\ee
It is easy to see that $K$ is compact,
the restriction  $T$ of $P_0\hgamma P_0-K$ to $\clh_{0s}$ is independent
of $s$, and $\sigma(T)=\sigma(\hgamma)$. Hence
$\sigma(P_0\hgamma P_0-K)=\sigma_{ess}(P_0\hgamma P_0-K)=\sigma(\hgamma)$.
Since $\sigma_{ess}(P_0\hgamma P_0-K)=\sigma_{ess}(P_0\hgamma P_0)$,
the proof follows.
\qed

\blmma\label{kernelgamma}
The operator $\gamma$ has trivial kernel.
In particular, for all  $r,s\in\half\bbz$,
$\ker \gamma_{rs}=\{0\}$.
\elmma
\prf
Let $P$ be the projection onto $\ker\gamma$.
Denote by $\phi$ the functional
$\phi:T\mapsto \langle e^{(0)}_{00}, T e^{(0)}_{00}\rangle$
Notice that the restriction of $\phi$ to $\cla$ is
the Haar state $h$, and therefore from appendix A1 in \cite{wo1}, we know
that 
\[
\phi(\chi_{\{q^{2n}\}}(\gamma))=(1-q^2)q^{2n}.
\]
Observe that $P=\chi_{\{0\}}(x)$.
Let $f(x)=\sum_{k=0}^n \chi_{\{q^{2k}\}}(x) +P $.
Then 
$\phi(f(\gamma))=1-q^{2n+2}+\phi(P)$.
Since
 $0\leq f(\gamma)\leq 1$,
we have $1-q^{2n+2}+\phi(P)\leq 1$ for all $n\in\bbn$.
Therefore $\phi(P)=0$, so that $Pe^{(0)}_{00}=0$.

From the commutation relations (\ref{suq2comm}), it follows that
$P\in \pi(\cla)'$. Since the vector $e^{(0)}_{00}$ is cyclic
for $\pi(\cla)$, it is separating for $\pi(\cla)'$. Therefore
we have $P=0$.
\qed

\bppsn\label{natureofgamma}
Let $r\in\half\bbn$ and $s\in\half\bbz$. The operator  $\gamma_{rs}$  is compact 
and its spectrum $\sigma(\gamma_{rs})$ coincides with
$\sigma(\gamma)=\{q^{2k}:k\in\bbn\}\cup\{0\}$.
Moreover  each $q^{2n}\in\sigma(\gamma_{rs})$ is
an eigenvalue of multiplicity 1.
\eppsn
\prf
Observe that
$\gamma_{rs}=P_{rs}(\gamma-\hgamma)P_{rs}+P_{rs}\hgamma P_{rs}$.
Since $\gamma-\hgamma$ is compact, and by lemma~\ref{hgamma}
$P_{rs}\hgamma P_{rs}$ is also compact, it follows that $\gamma_{rs}$
is compact.

Next we claim that $\sigma(\gamma_{rs})$ is independent of $s$.
 Write $T$ for the operator $P_{r,s+\half}\beta P_{rs}$. Then from the
commutation relations (\ref{alphacomm}--\ref{gmcomm}), it follows that $|T|=\gamma_{rs}^\half$
and the partial isometry $V_T$ appearing in the polar decomposition of
$T$ has initial space $\clh_{rs}$ and final space $\clh_{r,s+\half}$,
so that it can viewed as a unitary from $\clh_{rs}$ to $\clh_{r,s+\half}$.
Again from the commutation relations (\ref{alphacomm}--\ref{gmcomm}), we have
\bean
\gamma_{r,s+\half}T &=& (P_{r,s+\half}\gamma P_{r,s+\half})(P_{r,s+\half}\beta P_{rs})\\
&=& P_{r,s+\half}\gamma\beta P_{rs} \\
&=&P_{r,s+\half}\beta \gamma P_{rs} \\
&=& (P_{r,s+\half}\beta P_{rs})(P_{r,s+\half}\gamma P_{r,s+\half}) \\
&=& T \gamma_{r,s+\half}.
\eean
Therefore
\[
\gamma_{r,s+\half} V_T \gamma_{rs}^\half=V_T \gamma_{rs}^{\frac{3}{2}}.
\]
Since range of $\gamma_{rs}^\half$ is dense in $\clh_{rs}$, it
follows that $V_T^*\gamma_{r,s+\half} V_T =\gamma_{rs}$.
Thus $\gamma_{r,s+\half}$ and $\gamma_{rs}$ are unitarily equivalent.
So their spectrums are the same.

Since $P_r\gamma P_r=\oplus_s \gamma_{rs}$, it follows that
\be\label{spec1}
\sigma(P_r\gamma P_r)=\sigma_{ess}(P_r\gamma P_r)=\sigma(\gamma_{rs}).
\ee

Our next claim is that $\sigma_{ess}(P_r\gamma P_r)=\sigma(\gamma)$.
Let $K$ be the operator in the proof of lemma~\ref{hgamma}.
We have seen that $K$ is compact, 
$P_{0s}(P_0\hgamma P_0-K)P_{0s}$ is independent of $s$ and 
\[
\sigma(P_{0s}(P_0\hgamma P_0-K)P_{0s})=\sigma(\hgamma).
\]
Let $K_r:=U_r^* K U_r$. Then
\bean
P_{rs}(P_r\hgamma P_r-K_r)P_{rs} &=&
  U_r^*(U_rP_{rs}U_r^*(U_r P_r U_r^* U_r\hgamma U_r^* U_r P_r U_r^*-K)U_rP_{rs}U_r^*)U_r\\
&=& U_r^*(P_{0s}(P_0\hgamma P_0-K)P_{0s})U_r.
\eean
Therefore $P_{rs}(P_r\hgamma P_r-K_r)P_{rs}$ is independent of $s$ and
\[
\sigma(P_{rs}(P_r\hgamma P_r-K_r)P_{rs})=\sigma(P_{0s}(P_0\hgamma P_0-K)P_{0s})=\sigma(\hgamma).
\]
Hence
\[
\sigma_{ess}(P_r\hgamma P_r)=\sigma_{ess}(P_r\hgamma P_r-K_r)=\sigma(\hgamma)=\sigma(\gamma).
\]
Finally,
\[
\sigma_{ess}(P_r\gamma P_r)=\sigma_{ess}(P_r(\gamma-\hgamma) P_r+P_r\hgamma P_r)
=\sigma_{ess}(P_r\hgamma P_r)=\sigma(\gamma).
\]

It follows from (\ref{beta}) that
\[
\gamma (e^{(n)}_{ij}) = 
       k_{-1}(n,i,j) e^{(n-1)}_{ij}+ k_0(n,i,j) e^{(n)}_{ij}
               +k_1(n,i,j) e^{(n+1)}_{ij}, 
\]
where
\bean
k_1(n,i,j)  & =& -\Bigl(q^{4n+2i+2j+2}\frac{(1-q^{2n+2j+2})(1-q^{2n-2i+2})(1-q^{2n-2j+2})(1-q^{2n+2i+2})}
{(1-q^{4n+2})(1-q^{4n+4})(1-q^{4n+4})(1-q^{4n+6})}\Bigr)^\half,\cr
   && \label{kplus}\\
k_0(n,i,j)  & =& q^{2(n+j)}\frac{(1-q^{2n-2j})(1-q^{2n+2i})}{(1-q^{4n})(1-q^{4n+2})}
           +q^{2(n+i)}\frac{(1-q^{2n+2j+2})(1-q^{2n-2i+2})}{(1-q^{4n+2})(1-q^{4n+4})},\cr
   && \label{kzero}\\
k_{-1}(n,i,j)  & =& -\Bigl(q^{4n+2i+2j-2}\frac{(1-q^{2n-2j})(1-q^{2n+2i})
           (1-q^{2n+2j})(1-q^{2n-2i})}{(1-q^{4n-2})(1-q^{4n})(1-q^{4n})(1-q^{4n+2})}\Bigr)^\half.\cr
   && \label{kminus}
\eean
Therefore 
the operator $\gamma_{rs}-q^{2n}$ is a tridiagonal operator 
of the form
\[
 e^{(|r|+|s|+k)}_{s-r,-s-r}\mapsto 
  \cases{ b_0e^{(|r|+|s|)}_{s-r,-s-r}+c_0 e^{(|r|+|s|+1)}_{s-r,-s-r}  & if $k=0$,\cr
     a_k e^{(|r|+|s|+k-1)}_{s-r,-s-r}+ b_k e^{(|r|+|s|+k)}_{s-r,-s-r}+c_k e^{(|r|+|s|+k+1)}_{s-r,-s-r}  & if $k>0$,}
\]
with all the coefficients $a_k$, $b_k$ and $c_k$ nonzero.
It follows from this that the kernel of  $\gamma_{rs}-q^{2n}$
can have dimension at most 1.
Since $\gamma_{rs}$ is compact, each $q^{2n}\in\sigma(\gamma_{rs})$
is an eigenvalue. Therefore each $q^{2n}$ is an eigenvalue of multiplicity 1.
\qed

\newsection{The modular conjugation}
We will compute the modular conjugation operator for the Haar state 
in this section.
\bppsn\label{formula_S}
Denote by $S$ the operator $a\mapsto a^*$ on $\cla$.
Then viewed as an operator on $L_2(h)$, the set
$\{e^{(n)}_{ij}: (n,i,j)\in\Lambda\}$
is contained in the domain of $S$ and
\be\label{exprS}
Se^{(n)}_{ij}=(-1)^{2n+i+j}q^{i+j}e^{(n)}_{-i,-j}.
\ee
\eppsn
\prf
Recall (equation~57, page 115, \cite{k-s}) that if $t^{(n)}_{ij}$
denotes the ij\raisebox{.4ex}{th} matrix entry of the irreducible
representation indexed by $n$, then $e^{(n)}_{ij}$'s are
just the normalized  $t^{(n)}_{ij}$'s, more specifically,
\be\label{or}
e^{(n)}_{ij}=q^{-n+i} \Bigl(\frac{1-q^{4n+2}}{1-q^2}\Bigr)^{1/2} t^{(n)}_{ij}.
\ee
Therefore $\{e^{(n)}_{i,j}: (n,i,j)\in\Lambda\}$
is contained in the domain of $S$ and
\bean
S  e^{(n)}_{i,j}
&=&
\sum_{m,k,l}\langle e^{(m)}_{k,l}, S e^{(n)}_{i,j}\rangle e^{(m)}_{k,l}\\
&=&
\sum_{m,k,l}\langle e^{(n)}_{i,j}e^{(m)}_{k,l}, 1\rangle e^{(m)}_{k,l}.
\eean

By properties of Clebsch-Gordon coefficients,
\[
\langle e^{(n)}_{i,j}e^{(m)}_{k,l}, 1\rangle=0 \mbox{ for }n\neq m.
\]
From the equation preceeding equation~(43), p--74, \cite{k-s}, we get
\[
\sum_{l=0}^{2n} C_q(n,n,m;b-l+n,l-n,b)t^{(n)}_{a-k,b-l+n}t^{(n)}_{k,l-n}
=
C_q(n,n,m;a-k,k,a)t^{(m)}_{ab},\quad m=0,1,\ldots,2n.
\]
Here $C_q(m,n,p;i,j,k)$ are the Clebsch-Gordon coefficient 
$\langle e^{(p)}_k, e^{(m)}_i\otimes e^{(n)}_j\rangle$.
If we write
\[
A^{(b)}_{mj}=C_q(n,n,m; b-j+n,j-n,b),\quad m,j=0,1,\ldots,2n,
\]
then the above says that
\begin{displaymath}
\sum_{j=0}^{2n} A^{(b)}_{mj}t^{(n)}_{a-k,b-j+n}t^{(n)}_{k,j-n}=C_q(n,n,m;a-k,k,a)t^{(m)}_{ab}.
\end{displaymath}
for $m=0,1,\ldots,2n$.
Therefore
\[
t^{(n)}_{a-k,b-j+n}t^{(n)}_{k,j-n}
=
\sum_m ({A^{(b)}}^{-1})_{jm}C_q(n,n,m;a-k,k,a)t^{(m)}_{ab}.
\]
For $a=b=0$, coefficient of $t^{(0)}_{00}$ on the right hand side
is
\[
({A^{(0)}}^{-1})_{j0}C_q(n,n,0;-k,k,0),
\]
where $A^{(0)}$ is the matrix $(\!(A^{(0)}_{ij})\!)_{ij}$.
Thus,
\[
\langle t^{(n)}_{-k,-j+n}t^{(n)}_{k,j-n},1\rangle=
({A^{(0)}}^{-1})_{j0}C_q(n,n,0;-k,k,0).
\]
From equation~(\ref{or}),
\begin{displaymath}
\|t^{(n)}_{ij}\|=q^{-i}(2n+1)_q^{-\half}=
q^{n-i}\Bigl(\frac{1-q^{4n+2}}{1-q^2}\Bigr)^{-\half}.
\end{displaymath}
Hence
\begin{equation}\label{ip}
\langle e^{(n)}_{-k,-j+n}e^{(n)}_{k,j-n},1\rangle=
q^{-2n}\Bigl(\frac{1-q^{4n+2}}{1-q^2}\Bigr)({A^{(0)}}^{-1})_{j0}C_q(n,n,0;-k,k,0).
\end{equation}

Let us next find $({A^{(0)}}^{-1})_{j0}$.
Using equation~(73), page~81, \cite{k-s}, we get
\[
A^{(0)}_{mj}=C_q(n,n,m; -j+n,j-n,0)=(-1)^{2n+m}C_{q^{-1}}(n,n,m; j-n,-j+n,0).
\]
Write $D=(\!(d_{ij})\!)$ and $B=(\!(B_{ij})\!)$, where
$d_{ij}=\delta_{ij}(-1)^{2n+i}$ and
$B_{m,j}=C_{q^{-1}}(n,n,m; j-n,-j+n,0)$.
Then $A^{(0)}=DB$, so that
${A^{(0)}}^{-1}=B^{-1}D$.
Since $B$ has real entries, it follows from equations~(46) and (47), page~75,  \cite{k-s} that
$B$ is orthogonal.
Therefore ${A^{(0)}}^{-1}=B^tD$.
Hence
\[
({A^{(0)}}^{-1})_{j0}=(B^tD)_{j0}=(B^t)_{j0}d_{00}=B_{0j}(-1)^{2n}
=(-1)^{2n}C_{q^{-1}}(n,n,0; j-n,-j+n,0)
\]
Using equation~(73), page~81, \cite{k-s} and the equation preceeding equation~(68), page~81, \cite{k-s}, we get
\[
({A^{(0)}}^{-1})_{j0}=C_q(n,n,0; -j+n,j-n,0)
=(-1)^jq^{2n-j}\Bigl(\frac{1-q^{4n+2}}{1-q^2}\Bigr)^{-\half},
\]
and
\[
C_q(n,n,0; -k,k,0)=(-1)^{n+k}q^{n-k}\Bigl(\frac{1-q^{4n+2}}{1-q^2}\Bigr)^{-\half}.
\]

Substituting these values in (\ref{ip}), we  get
\begin{equation}
\langle e^{(n)}_{ij}e^{(n)}_{-i,-j},1\rangle=
(-1)^{2n+i+j}q^{i+j}.
\end{equation}

Thus we have equation~(\ref{exprS}).
\qed

In particular, it follows from the above proposition that
the operator $S$ is closable. Let $\bar{S}$ denote the closure
of $S$.
Let $J$ denote the antilinear operator, given on the basis elements by
\begin{equation}\label{modular}
Je^{(n)}_{ij}=(-1)^{2n+i+j}e^{(n)}_{-i,-j},
\end{equation}
and let $\Delta$ be given by
\be
\Delta e^{(n)}_{ij}=q^{2i+2j}e^{(n)}_{ij}.
\ee
Then it follows from (\ref{exprS}) that $\bar{S}=J\Delta^\frac{1}{2}$,
and $\mbox{Dom}\,  \bar{S}= \mbox{Dom}\, \Delta^\half$. 
By lemma~1.5, \cite{t}, it follows that
$J$ is the modular conjugation and
$\Delta$ is the modular operator associated with the Haar state,
and by theorem~1.19, \cite{t}, we have
$\pi(\cla)'=J\pi(\cla)''J$.

\newsection{The main theorem}

Let $D$ be the operator given by (\ref{genericd}),
and let $F=\mbox{sign}\,D$. 
\bthm\label{main}
Let $T\in \pi(\cla)'$.
If $[F, T]$ is compact, then $T$ is a scalar.
\ethm
Let $J$ be the modular conjugation operator computed in the previous section.
Then the theorem says that if $T\in \pi(\cla)''$ and $[F, JTJ]$ is compact, then
$T$ must be a scalar.
We will first prove the following special case of the above theorem.

\bthm\label{aux}
Let $f$ be a complex valued function on $\sigma(\gamma)$.
If $[F, Jf(\gamma)J]$ is compact, then $f(\gamma)$ must be a scalar.
\ethm
We will need the following simple lemma for the proof of this
theorem.

\blmma\label{distance}
Let $A$ and $B$ be two compact operators with trivial kernel
such that $\sigma(A)=\sigma(B)$ and each nonzero element of
$\sigma(A)$ is an eigenvalue of multiplicity 1 for both
$A$ and $B$. 
Let $u$ be a unit eigenvector of $A$ corresponding
to an eigenvector $\lambda$, $v$ be a unit eigenvector of $B$
corresponding to the same eigenvalue $\lambda$.
Assume that $|\langle u,v\rangle|<1-\epsilon$ where $\epsilon>0$.
Then there is a positive constant $c=c(\epsilon,\lambda,\sigma(A))$
such that 
\[
 \|A-B\|\geq c.
\]
\elmma
\prf
It follows from the given conditions that there is a unitary
$U$ such that $Uu=v$ and $B=UAU^*$.
Let $w$ be the projection of $v$ onto $u^\perp$, i.e.
\[
 w=v-\langle u,v\rangle u.
\]
Then
\[
\|(A-B)v\|=\|(A-\lambda)v\|=\|(A-\lambda)w\|.
\]
Since $(A-\lambda)$ is invertible on $u^\perp$, it follows that
\[
\|w\|=\|((A-\lambda)|_{u^\perp})^{-1}(A-\lambda)w\|\leq
         \|((A-\lambda)|_{u^\perp})^{-1}\|\|(A-\lambda)w\|,
\]
so that 
\[
\|A-B\|\geq \|(A-B)v\|=\|(A-\lambda)w\| \geq \|w\| \|((A-\lambda)|_{u^\perp})^{-1}\|^{-1}.
\]
Observe that
\[ 
\|((A-\lambda)|_{u^\perp})^{-1}\|=
  (\inf\{|\lambda-\mu|:\mu\in\sigma(A), \mu\neq\lambda\})^{-1}<\infty.
\]
Since
\[
1=\|v\|^2=|\langle u,v\rangle|^2+\|w\|^2< (1-\epsilon)^2+\|w\|^2,
\]
the result follows.
\qed

\noindent \textit{Proof of theorem~\ref{aux}}:
Write $T:=f(\gamma)$.
Let $Q=\frac{I-JFJ}{2}$. Then compactness of $[F,JTJ]$ is equivalent
to compactness of $[Q,T]$.
Let  $w_{rs}=e^{|r|+|s|}_{s-r,-s-r}$.
For an operator $A$ on $\clh$, denote by $A_{rs}$ the operator $P_{rs}A P_{rs}$.
The projection $Q$ commutes with each $P_{rs}$ and one has
\[
 Q_{rs}=\cases{|w_{rs}\rangle\langle w_{rs}| & if $r, -s\in\half\bbn$,\cr
      0 & otherwise.}
\]
Since $\gamma$ also commutes with $P_{rs}$, we have
\[
  [Q,T]=\oplus_{r,s}[Q_{rs},T_{rs}]=\oplus_{r,-s\in\half\mathbb{N}}[Q_{rs},T_{rs}].
\]
Take $\ell\neq m\in\bbn$ and let $\lambda=f(q^{2\ell})$, $\mu=f(q^{2m})$.
Recall that for each $r\in\half\bbn$, $s\in\half\bbz$,
one has $\sigma(\gamma_{rs})=\{q^{2n}:n\in\bbn\}\cup\{0\}$,
and each $q^{2n}$ is an eigenvalue of multiplicity 1.
Let $u_{rs}$ and $v_{rs}$ be eigenvectors of $\gamma$ in $\mathcal{H}_{rs}$
corresponding to the eigenvalues $q^{2\ell}$ and $q^{2m}$ respectively.
Then $Tu_{rs}=\lambda u_{rs}$ and $Tv_{rs}=\mu v_{rs}$.
Now 
\bean
\langle u_{rs}, [Q,T](u_{rs}-v_{rs}) \rangle &=&
   \langle u_{rs},Q(\lambda u_{rs}-\mu v_{rs})\rangle-\langle T^* u_{rs},Q(u_{rs}-v_{rs})\rangle\\
&=&     \lambda \langle u_{rs},Q u_{rs}\rangle -
  \mu \langle u_{rs}, Q v_{rs}\rangle - \lambda \langle u_{rs}, Q u_{rs}\rangle
  + \lambda \langle u_{rs}, Q v_{rs}\rangle  \\
&=& (\lambda-\mu)\langle u_{rs}, Q v_{rs}\rangle.
\eean
In particular,
\bea
\|[Q,T](u_{r0}-v_{r0})\|^2  &\geq & |\lambda-\mu|^2 |\langle u_{r0}, Q v_{r0}\rangle|^2\\
&=& |\lambda-\mu|^2 
   |\langle u_{r0}, w_{r0}\rangle   \langle w_{r0}, v_{r0}\rangle|.\label{ineq1}
\eea
We will show that for some subsequence $r_k$,
\bea
\lim_{k\rightarrow\infty}\langle u_{r_k 0}, w_{r_k 0}\rangle &\neq & 0,\label{nonzero1}\\
\lim_{k\rightarrow\infty}\langle v_{r_k 0}, w_{r_k 0}\rangle &\neq & 0.\label{nonzero2}
\eea
Observe that 
\[
 \langle v_{r0},w_{r0}\rangle = \langle U_r v_{r0}, U_r w_{r0}\rangle
    = \langle U_r v_{r0}, w_{00}\rangle.
\]
Now $U_r v_{00}$ be an eigenvector of $U_r \gamma_{r0}U^*$ corresponding
to the eigenvalue $q^{2m}$.
Let $\xi$ be a unit eigenvector  of $\hgamma$ in $\clh_{00}$ corresponding
to the eigenvalue $q^{2m}$.
Note that
\[
 U_r \gamma_{r0} U_r^*-\hgamma_{00} 
   = U_r(\gamma_{r0}-\hgamma_{r0})U_r^* + U_r\hgamma_{r0} U_r^*-\hgamma_{00} 
  = U_r(\gamma_{r0}-\hgamma_{r0})U_r^*.
\]
Therefore
\[
 \|U_r \gamma_{r0} U_r^*-\hgamma_{00}\|=\|\gamma_{r0}-\hgamma_{r0}\|
   =\|P_{r0}(\gamma-\hgamma)P_{r0}\|.
\]
Since $\gamma-\hgamma$ is compact, it follows that
\[
\lim \|U_r \gamma_{r0} U_r^*-\hgamma_{00}\|=0.
\]
By lemma~\ref{distance}, 
\[
 \lim|\langle U_r v_{r0},\xi\rangle|=1.
\]
Let $\theta_r\in[0,2\pi]$ be such that
$\exp(i\theta_r)=\langle U_r v_{r0},\xi\rangle$.
Write $\xi_r=\exp(-i\theta_r)\xi$.
Then
$\langle U_r v_{r0},\xi_r\rangle=1$ for each $r$.
Choose a subsequence $r_k$ such that $\theta_{r_k}$
converges, to say $\theta$. Let $\zeta=\exp(-i\theta)\xi$.
Then it follows that
\[
 \lim_k \langle U_{r_k}v_{r_k 0},\zeta\rangle =1.
\]
Therefore 
\[
 \lim_k \|U_{r_k}v_{r_k 0}-\zeta\|=0.
\]
Hence
\[
\langle v_{r_k 0}, w_{r_k 0}\rangle =\langle U_{r_k}v_{r_k 0}, w_{00}\rangle
  = \langle U_{r_k}v_{r_k 0}-\zeta , w_{00} \rangle +\langle\zeta, w_{00}\rangle.
\]
The first term converges to zero. Let us show that the second term is nonzero.
Let 
\[
k=\min\{n\in\bbn: \langle \zeta, e^{(n)}_{00}\rangle\neq 0\}.
\]
Then one has 
\[
 \langle \zeta,\hgamma^k w_{00}\rangle = \langle \zeta,\hgamma^k e^{(0)}_{00}\rangle\neq 0.
\]
Since $\hgamma\zeta=q^{2m}\zeta$,
we have
\[
\langle  \zeta, w_{00}\rangle = 
   q^{-2mk}\langle \hgamma^k \zeta, w_{00}\rangle
 = q^{-2mk}\langle  \zeta,\hgamma^k w_{00}\rangle\neq 0.
\]
Thus we have (\ref{nonzero2}).
An identical proof shows that $r_k$ will have a further subsequence,
which we continue to denote by $r_k$ by abuse of notation,
for which we have both (\ref{nonzero2}) and (\ref{nonzero1}).

Since $[Q,T]$ is compact, $[Q,T](u_{r_k 0}-v_{r_k 0})$ converges to zero.
Therefore by (\ref{ineq1}), 
we must have $\lambda=\mu$,
i.e.\ $f(q^{2\ell})=f(q^{2m})$.
Since this is true for all $\ell \neq m \in\bbn$,
$T=f(\gamma)$ must be a scalar.
\qed

\bppsn\label{aux2}
Let $m,n\in\bbz$ and let $f$ be a complex valued function on $\sigma(\gamma)$.
Assume $m \neq 0$.
If $[F, J\alpha_m \beta_n f(\gamma)J]$ is compact, then 
$\alpha_m \beta_n f(\gamma)=0$.
\eppsn
\prf
Assume $m>0$, so that $\alpha_m=\alpha^m$.
Compactness of $[F, J\alpha_m \beta_n f(\gamma)J]$ implies
compactness of $[Q, \alpha_m \beta_n f(\gamma)]$.
Since $Q$ is self-adjoint,  this implies
$[Q,  (f(\gamma))^*\beta_n^*(\alpha^m)^* \alpha^m \beta_n f(\gamma))]$ is compact.
Now $\beta_n^*(\alpha^m)^* \alpha^m \beta_n$ is of the form
$p(\gamma)$ for some polynomial $p$.
By theorem~\ref{aux} it follows that 
$(f(\gamma))^*\beta_n^*(\alpha^m)^* \alpha^m \beta_n f(\gamma))$
is a scalar. Suppose it is nonzero. Since $(\alpha^m)^* \alpha^m$
is a polynomial in $\gamma$, we have 
\[
(f(\gamma))^*\beta_n^*(\alpha^m)^* \alpha^m \beta_n f(\gamma))
   =(f(\gamma))^*\beta_n^*\beta_n f(\gamma)) (\alpha^m)^* \alpha^m
\]
It would then follow that
the kernel of $(\alpha^m)^* \alpha^m$ is trivial. 
This implies that the kernel of $\alpha$ is trivial. But
$\ker \alpha=\ker \alpha^*\alpha = \ker (1-\gamma)$ which
is infinite dimensional by proposition~\ref{natureofgamma}.
Therefore we must have
$(f(\gamma))^*\beta_n^*\alpha_m^* \alpha_m \beta_n f(\gamma))=0$
which implies $\alpha_m \beta_n f(\gamma)=0$.

For $m<0$, observe that $(\alpha_m \beta_n f(\gamma))^*=\alpha_{-m}\beta_{-n}g(\gamma)$
for some function $g$ and use the above argument.
\qed

\bppsn\label{aux3}
Let $n\in\bbz$ and let $f$ be a nonzero complex valued function on $\sigma(\gamma)$.
If $n \neq 0$, then $[F, J \beta_n f(\gamma)J]$ is not compact.
\eppsn
\prf
As before, compactness of $[F, J \beta_n f(\gamma)J]$  is equivalent to the compactness of
$[Q, \beta_n f(\gamma)]$. So it is enough to show that
$[Q, \beta_n f(\gamma)]$ is not compact.
Also it is enough to prove this for $n>0$.

Since $\beta$ and $\beta^*$ both have trivial kernel,
the partial isometry $V$ appearing in the polar decomposition
of $\beta$ is unitary and
$\beta_n f(\gamma)=V^n g(\gamma)$ for some function $g$.
Let $m\in\bbn$ be such that $\lambda:=g(q^{2m})\neq 0$.
For $r\in\half\bbn$, let $v_{r0}$ be a unit eigenvector
of $\gamma$ in $\clh_{r0}$ corresponding to the eigenvalue $q^{2m}$.
For $s\in\half\bbz$, define 
\[
v_{rs}:= V^{2s}v_{r0}.
\]
Then  $v_{rs}$ is a unit vector in $\clh_{rs}$ and
since $V$ commutes with $\gamma$, we have $\gamma v_{rs}=q^{2m} v_{rs}$.
Therefore $g(\gamma)v_{rs}=\lambda v_{rs}$ for all $r,s$.
We then have
\bean
\langle v_{r,\frac{n}{2}} , [Q, V^n g(\gamma)]v_{r0}\rangle
  &=& \langle v_{r,\frac{n}{2}} , QV^n g(\gamma)v_{r0}\rangle
         - \langle v_{r,\frac{n}{2}} , V^n g(\gamma) Q v_{r0}\rangle  \\
&=& \lambda \langle v_{r,\frac{n}{2}} , QV^n v_{r0}\rangle
          - \langle g(\gamma)v_{r,\frac{n}{2}} , V^n  Q v_{r0}\rangle  \\
&=& \lambda (\langle v_{r,\frac{n}{2}} , Qv_{r,\frac{n}{2}}\rangle
               -  \langle v_{r0},  Q v_{r0}\rangle )\\
&=& -\lambda\|Q v_{r0}\|^2\\
&=& -\lambda |\langle v_{r0}, w_{r0}\rangle|^2.
\eean
From the proof of theorem~\ref{aux}, there is a sequence $r_k$ such that
$\lim_k \langle v_{r_k 0}, w_{r_k 0}\rangle\neq 0$. Therefore
the operator $[Q, V^n g(\gamma)]$ can not be compact.
\qed

We now have all the ingredients ready for the proof of theorem~\ref{main}.
In order to make use of these, 
we need to look at certain operator valued Fourier coefficients.

Let $\tau$ be the action of $S^1\times S^1$  on $\cla$ by automorphisms given by
\begin{displaymath}
\tau_{z,w}:\cases{\alpha\mapsto z\alpha,\cr
                    \beta\mapsto w\beta.}
                \end{displaymath}
Let $V_{z,w}:L_2(h)\rightarrow L_2(h)$ be given by
\[
V_{z,w}e^{(n)}_{ij}=z^{-i-j}w^{i-j}e^{(n)}_{ij}.
\]
Then $\pi(\tau_{z,w}(a))=V_{z,w}\pi(a) V_{z,w}^*$ for all $a\in\cla$.
Thus the action extends to a strongly continuous action of $S^1\times S^1$ on
the von Neumann algebra $\pi(\cla)''$.
For $T\in \pi(\cla)''$ and $m,n\in\bbz$, denote by $\scrf_{mn}(T)$ the following operator:
\begin{displaymath}
\scrf_{mn}(a)=\int_{S^1}\int_{S^1}z^{-m}w^{-n}\tau_{z,w}(T)dz\,dw.
\end{displaymath}
Note that the above integral is defined in the strong sense. 
In case the integrand is norm continuous, it
coincides with the corresponding integral in the norm sense.
\blmma\label{convfourier}
Let $T_k$ be a sequence of operators in $\pi(\cla)''$ that converges strongly to
an operator $T$. Then for all $m,n\in\bbz$, the sequence
$\scrf_{mn}(T_k)$ converges strongly to $\scrf_{mn}(T)$.
\elmma
\prf
Take a vector $u\in L_2(h)$.
Then
\[
 \scrf_{mn}(T_k)u=\int\int z^{-m}w^{-n}V_{z,w}T_k V_{z,w}^*u\, dz\, dw, \quad
\scrf_{mn}(T)u=\int\int z^{-m}w^{-n}V_{z,w}T V_{z,w}^*u \,dz\, dw.
\]
Since $T_k$ converges strongly to $T$, for each $z,w\in S^1$, we have
\[
 \lim_{k\rightarrow\infty}z^{-m}w^{-n}V_{z,w}T_k V_{z,w}^*u
          = z^{-m}w^{-n}V_{z,w}T V_{z,w}^*u,
\]
and
\bean
\|z^{-m}w^{-n}V_{z,w}T_k V_{z,w}^*u\| &\leq &
       \|T_k\|\|u\|  \\
&\leq & (\sup_k \|T_k\|)\|u\|.
\eean
Now an application of the Dominated Convergence Theorem
for Banach space valued functions (Theorem~3, page--45, \cite{d-u})
gives us the required result.
\qed

\blmma\label{fourier}
Let $T\in \pi(\cla)''$. If $\scrf_{mn}(T)=0$ for all $(m,n)\neq (0,0)$,
then $T=f(\gamma)$ for some bounded measurable function
$f$ on $\sigma(\gamma)$.
\elmma
\prf
Let $B$ be the *-subalgebra of $\cla$ consisting of finite linear
combinations of elements of the form $\alpha_m\beta_n \gamma^k$,
where $m,n\in\bbz$ and $k\in\bbn$.
Clearly $B$ is dense in $\cla$.
Observe that 
\begin{enumerate}
\item 
for any $T\in B$,
one has $\scrf_{mn}(T)=\alpha_m\beta_n p(\gamma)$ for some
polynomial $p$,
\item
if $T=\alpha_m\beta_n p(\gamma)$ for some polynomial $p$, then
\[
\scrf_{jk}(T)=\cases{ T & if $j=m$, $k=n$,\cr
                  0 & otherwise.}
\]
\end{enumerate}
Now let $T\in \pi(\cla)''$ with $\scrf_{mn}(T)=0$ for all $(m,n)\neq (0,0)$.
Take any two vectors $u$ and $v$ in $\clh$ and let
$f:S^1\times S^1\rightarrow \bbc$ be the function given by
$f(z,w)=\langle u, \tau_{z,w}(T)v\rangle$.
Then by the above condition on $T$, it follows that
all the Fourier coefficients $\hat{f}(m,n)$ are zero for
all $(m,n)\neq (0,0)$. This implies $f$ is a constant function.
Since this is true for any two vectors $u$ and $v$,
it follows that $(z,w)\mapsto \tau_{z,w}(T)$ is constant,
so that for all $z,w\in S^1$, we have $\tau_{z,w}(T)=T$.
Therefore $\scrf_{00}(T)=T$.
Let $T_k$ be a sequence in $B$ that converges strongly to $T$.
By lemma~\ref{convfourier}, we have s-$\lim_k\scrf_{mn}(T_k)=\scrf_{mn}(T)$
for all $m,n\in\bbz$.
In particular, we have s-$\lim_k\scrf_{00}(T_k)=\scrf_{00}(T)=T$.
Since each $\scrf_{00}(T_k)$ is of the form $p_k(\gamma)$ for some
polynomial $p_k$, the operator $T$ must be of the form
$f(\gamma)$ for some bounded measurable function on $\sigma(\gamma)$.
\qed

\blmma
Let $T\in \pi(\cla)''$. Then for $m,n\in\bbz$,
the operator $\scrf_{mn}(T)$ is of the form
$\alpha_m\beta_n f(\gamma)$ for some  function
$f$ on $\sigma(\gamma)$.
\elmma
\prf
Since $(\scrf_{mn}(T))^*=\scrf_{-m,-n}(T^*)$, it is enough
to prove the statement for $m\leq 0$.
So assume $m\leq 0$.
Let $T_k$ be a sequence in $B$ that converges strongly to $T$.
Then by lemma~\ref{convfourier}, $\scrf_{mn}(T_k)$ converges
strongly to $\scrf_{mn}(T)$. Each $\scrf_{mn}(T_k)$ is of the form
$\alpha_m\beta_n p_k(\gamma) =(\alpha^*)^{|m|}\beta_n p_k(\gamma)$ for some polynomial $p_k$.
Now
\[
((\alpha^*)^{|m|}\beta_n)^*((\alpha^*)^{|m|}\beta_n)
   = \alpha^{|m|}(\alpha^*)^{|m|}\gamma^{|n|},
\]
and $\ker\alpha^*=\{0\}=\ker\beta=\ker\beta^*$.
Therefore the operator $(\alpha^*)^{|m|}\beta_n$
has trivial kernel. Therefore the polar decomposition of
$(\alpha^*)^{|m|}\beta_n$ is of the form $V\sqrt{r(\gamma)}$
where $V$ is an isometry and $r$ is a polynomial.
Thus $V\sqrt{r(\gamma)}p_k(\gamma)$ converges strongly to $\scrf_{mn}(T)$.
Therefore $\sqrt{r(\gamma)}p_k(\gamma)$ converges strongly to $V^*\scrf_{mn}(T)$.
It follows that $V^*\scrf_{mn}(T) = f(\gamma)$ for some bounded function $f$
and $\lim_k \sqrt{r(x)}p_k(x)=f(x)$ for all $x\in\sigma(\gamma)$.
Define functions $\widetilde{p}$ and $\widetilde{p}_k$ on $\sigma(\gamma)$
as follows:
\[
\widetilde{p}(x)=\cases{f(x)/\sqrt{r(x)} & if $r(x)\neq 0$,\cr
        0 & if $r(x)=0$,} \qquad
\widetilde{p}_k(x)=\cases{p_k(x) & if  $r(x)\neq 0$,\cr
        0 & if $r(x)=0$.}
\]
Then $\sqrt{r(x)}p_k(x)=\sqrt{r(x)}\widetilde{p}_k(x)$
and
$f(x)=\sqrt{r(x)}\widetilde{p}(x)$.
This means $\sqrt{r(\gamma)}p_k(\gamma)=\sqrt{r(\gamma)}\widetilde{p}_k(\gamma)$
and $\sqrt{r(\gamma)}p_k(\gamma)$ converges strongly to
$\sqrt{r(\gamma)}\widetilde{p}(\gamma)$.
Therefore
$V\sqrt{r(\gamma)}p_k(\gamma)$ converges strongly to
$V\sqrt{r(\gamma)}\widetilde{p}(\gamma)=(\alpha^*)^{|m|}\beta_n\widetilde{p}(\gamma)$.
Hence $\scrf_{mn}(T)=(\alpha^*)^{|m|}\beta_n\widetilde{p}(\gamma)$.
\qed

\noindent We now turn to the proof of theorem~\ref{main}.\\[1ex]
\textit{Proof of theorem~\ref{main}}:
Compactness of $[F, JTJ]$ implies
$[Q, T]$ is compact. Since $V_{z,w}[Q,T]V_{z,w}^*=[Q,\tau_{z,w}(T)]$,
it follows that $[Q,\scrf_{mn}(T)]$ is compact for all $m$ and $n$.
Since the operator $\scrf_{mn}(T)$ is of the form
$\alpha_m\beta_n f(\gamma)$ for some  function
$f$ on $\sigma(\gamma)$,
by propositions~\ref{aux2} and \ref{aux3} we get $\scrf_{mn}(T)=0$ for all $(m,n)\neq (0,0)$.
An application of lemma~\ref{fourier} now tells us that $T=f(\gamma)$ for some
bounded function $f$. Hence using theorem~\ref{aux}, we get that
$T$ is a scalar.
\qed

\brmrk\rm
By the characterization of equivariant spectral triples in \cite{c-p}
(see the discussion preceeding proposition~4.4, \cite{c-p}),
for any equivariant $D$, $\mbox{sign}\,D$ has to be of the form
$2P-I$ or $I-2P$, where $P$ is the projection onto the subspace
spanned by
$\{e^{(n)}_{ij}: n\in\half\bbn, n-i\in E, j=-n,-n+1,\ldots,n\}$,
$E$ being some finite subset of $\bbn$. A slight modification in the
proof of theorem~\ref{aux} will work for the sign of any such $D$.
\ermrk

\bcrlre\label{main2}
Suppose $T\in \pi(\cla)'$. If $[D,T]$ is bounded, then $T$ must be a scalar.
\ecrlre
\prf
Boundedness of $[D,T]$ implies compactness of $[F,T]$. Therefore
the result follows from theorem~\ref{main}.
\qed

\brmrk\rm
Suppose $\pi$ is a faithful representation of a $C^*$-algebra $\cla$ on a Hilbert space
$\mathcal{H}$ and $(\mathcal{H},\pi,D)$ is a spectral triple for $\cla$.
If there is another $C^*$-algebra $\clb$ and a faithful representation $\rho$ of $\clb$ on $\mathcal{H}$
such that $\pi(a)$ and $\rho(b)$ commute  for all $a\in \cla$, $b\in \clb$
and $(\mathcal{H},\rho,D)$ is a spectral triple for $\clb$, then the pair $(\mathcal{H},D)$ together with the representation $\pi\otimes\rho: a\otimes b\mapsto \pi(a)\otimes\rho(b)$
gives rise to a spectral triple for $\cla\otimes \clb$ and hence an element
in the $K$-homology of $\cla\otimes \clb$.

What the above corollary says is that for the Dirac operator constructed
in \cite{c-p} on $L_2(h)$ along with the representation by left
multiplication, such a pair $(\clb,\rho)$ does not exist (other than
the trivial one: $\clb=\bbc$), thereby
preventing one from turning it into a spectral triple for 
$\mathcal{A}\otimes \clb$ in a natural manner.
Thus the triple $(L_2(h), \pi, D)$ does not admit a rational 
Poincar\'e dual in the sense of Moscovici (\cite{m}).
\ermrk

\brmrk\label{rem412}\rm
Note that the Dirac operator we have considered here is the
one equivariant with respect to the right regular representation
of the group $SU_q(2)$. Recall (\cite{c-p4}) that a generic Dirac operator equivariant
with respect to the left regular representation is of the form
$D:e^{(n)}_{ij}\mapsto d(n,i)e^{(n)}_{ij}$, where
\[
 d(n,i)=\cases{2n+1 & if $n\neq j$,\cr
                      -(2n+1) & if $n=j$.}
\]
All the results in this section continue to hold for this Dirac
operator as well. The proofs also go through verbatim, except the proof of
proposition~\ref{aux3}, where one has to look at commutators
$[F, J\beta_n f(\gamma)J]$ for $n<0$.
\ermrk

\section{K-theory fundamental class}
We will show in this section that even though the spectral triple
we considered does not give a fundamental class, a little modification
enables one to construct a fundamental class that gives Poincar\'e duality.

We start the section with the following straightforward
but important observation.
\bthm\label{pd}
Poincar\'e duality holds for $\mathcal{A}$.
\ethm
\prf
It follows from the description of the 
irreducible representations of $\mathcal{A}$~(\cite{ko-so})
that it is a type I $C^*$-algebra.
Both $\mathcal{A}$ and $C(S^1)$ are separable
type~I  $C^*$-algebras
and have the same $K_0$ and $K_1$ groups.
Therefore it follows from Rosenberg \& Schochet (\cite{r-s} that
$\mathcal{A}$ and $C(S^1)$ are
$KK$-equivalent.
Poincar\'e duality holds for $C(S^1)$, hence it 
follows from lemma~3.4, \cite{b-m-r-s}
that Poincar\'e duality holds for $\mathcal{A}$  also.
\qed

One can see that Poincar\'e duality is just a consequence of the $KK$-theoretic
properties of the underlying $C^*$-algebra. What is of greater interest
is to get an explicit realization of the $K$-homology fundamental class. 
Thus we want to identify explicitly a class in
$KK(\mathcal{A}\otimes \mathcal{A},\bbc)$ that will give us a 
$K$-homology fundamental class.
As a first step, we will exhibit an element in $KK(\mathcal{A},C(S^1))$
that will give us a $KK$-equivalence. We then compose this with
the fundamental class for the torus to construct the desired class.
This involves computing the Kasparov product of two elements, which can
sometimes be difficult. As we will see, we avoid computing any nontrivial
Kasparov product by exploiting the special form of the $KK$-equivalence
we construct.

\blmma\label{tech1}
$KK(\mathcal{A},\bbc)\otimes KK(\bbc,C(S^1))\cong KK(\mathcal{A}, C(S^1))$.\\
(here $KK(A,B)$ means $KK_0(A,B)\oplus KK_1(A,B)$ )
\elmma
\prf
Observe that $KK(\mathcal{A},\bbc)=\bbz\oplus\bbz$ and 
$KK(\bbc, C(S^1))=\bbz\oplus\bbz$.
Thus both are torsion-free and by the K\"{u}nneth theorem (due to Rosenberg \& Schochet,
theorem~23.1.2, \cite{b}), the result follows.
\qed

\blmma\label{tech2}
$KK(\mathcal{A}, C(S^1))\cong M_2(\bbz)$.
\elmma
\prf
Since $KK(\bbc, \mathcal{A})=\bbz\oplus\bbz$, by the universal coefficient theorem (UCT) (theorem~23.1.1, \cite{b})
it follows that
\[
KK(\mathcal{A},C(S^1))\cong Hom(KK(\bbc, \mathcal{A}),KK(\bbc, C(S^1)))
\cong Hom (\bbz\oplus\bbz,\bbz\oplus\bbz)\cong M_2(\bbz).
\]
Note that in the above isomorphism, an element $\eta\in KK(\mathcal{A},C(S^1))$
is mapped to the homomorphism given by
$\xi\mapsto \xi\otimes\eta$.
\qed

\blmma
Let $\xi_i \in K_i(\mathcal{A})$, $\zeta_i\in K^i(\mathcal{A})$ and $\eta_i\in K_i(C(S^1))$,
$i=0,1$. Define
\[ 
\gamma:=\zeta_0\otimes\eta_0 + \zeta_1\otimes\eta_1.
\]
Then the map $\xi\mapsto \xi\otimes \gamma$
takes $\xi_0$ to $\langle \xi_0 ,\zeta_0\rangle  \eta_0$ and 
$\xi_1$ to $\langle \xi_1,\zeta_1\rangle  \eta_1$.
\elmma
\prf
Proof follows immediately from the observations that
$\xi_0\otimes\zeta_1\equiv\langle\xi_0,\zeta_1\rangle$
and 
$\xi_1\otimes\zeta_0\equiv\langle\xi_1,\zeta_0\rangle$
are both zero, being elements of $K_1(\bbc)$.
\qed

\bppsn
Let $\sigma$ denote the trivial grading on $\bbc$.
Then the even Fredholm module $(\bbc,\sigma,\epsilon,0)$ gives
a generator for $K^0(\mathcal{A})=\bbz$.
\eppsn
\prf
Since $\mathcal{A}$ and $C(S^1)$ are $KK$-equivalent, one has
$K^0(\mathcal{A})=\bbz$.
This together with the simple observation that
the pairing
$\langle [ (\bbc,\sigma,\epsilon,0)],[1]\rangle$
is 1
gives us the required result.
\qed

We now put together the two results above to
produce a $KK$-equivalence.

\bppsn\label{tech3}
Let $\zeta_1$ be the K-homology class of the equivariant triple
for $\mathcal{A}$ (under the $SU_q(2)$ action), i.e.\ 
$\zeta_1=[(L_2(h),\pi,D)]$.
Let $\eta_1$  be the element $[z]$ in $K_1(C(S^1))$.
Let $\zeta_0$ and $\eta_0$ be generators for 
$K^0(\mathcal{A})$ and $K_0(C(S^1))$ respectively.
Then 
$\gamma:=\zeta_0\otimes\eta_0 + \zeta_1\otimes\eta_1$
gives a $KK$-equivalence between $\mathcal{A}$ and $C(S^1)$.
\eppsn
\prf
Recall that $\gamma$ corresponds to the element $\xi\mapsto\xi\otimes\gamma$
in $Hom (KK(\bbc, \mathcal{A}),KK(\bbc, C(S^1)))$
and $KK(\bbc, \mathcal{A})\cong \bbz\oplus\bbz$, $KK(\bbc, C(S^1))\cong\bbz\oplus\bbz$.
Using these identifications, it is now easy to see that $\gamma$ maps the 
element $1\oplus 0$ to $1\oplus 0$
and $0\oplus 1$ to $0\oplus 1$. In other words, $\gamma$ is the element
$\pmatrix{1&0\cr 0&1}$ in $KK(\mathcal{A},C(S^1))\cong M_2(\bbz)$.
This being an invertible element, gives a $KK$-equivalence.
\qed

\bthm\label{fclass}
The spectral triple $(L_2(h)\oplus L_2(h), \pi\oplus\epsilon, D\oplus D)$  gives a 
fundamental class for $\mathcal{A}=C(SU_q(2))$.
\ethm
\prf
Let $\rho$ be the representation of $C(S^1)\otimes C(S^1)$
on $L_2(S^1)$ given by
\[
 \rho(f\otimes g) h= fgh,
\]
and let $\partial=\partial_\theta$ be the derivative.
Then $(L_2(S^1),\rho,\partial)$ gives the standard 
fundamental class for $C(S^1)$.

Write $\lambda$ for the class of $(L_2(S^1),\rho,\partial)$ in 
$KK^1(C(S^1)\otimes C(S^1),\bbc)$.
Then $(\gamma\otimes\gamma)\otimes\lambda$
gives a $K$-homology fundamental class for $\mathcal{A}$.
Now
\bean
\gamma\otimes\gamma &=& 
  (\zeta_0\otimes\eta_0+\zeta_1\otimes\eta_1)\otimes
      (\zeta_0\otimes\eta_0+\zeta_1\otimes\eta_1)\\
&=&  (\zeta_0\otimes\eta_0)\otimes(\zeta_0\otimes\eta_0)
 + (\zeta_0\otimes\eta_0)\otimes(\zeta_1\otimes\eta_1)\\
&& + (\zeta_1\otimes\eta_1)\otimes(\zeta_0\otimes\eta_0)
 + (\zeta_1\otimes\eta_1)\otimes(\zeta_1\otimes\eta_1)\\
&=& (\zeta_0\otimes\zeta_0)\otimes(\eta_0\otimes\eta_0)
 + (\zeta_0\otimes\zeta_1)\otimes(\eta_0\otimes\eta_1)\\
&& + (\zeta_1\otimes\zeta_0)\otimes(\eta_1\otimes\eta_0)
 + (\zeta_1\otimes\zeta_1)\otimes(\eta_1\otimes\eta_1).
\eean
Clearly $(\eta_0\otimes\eta_0)\otimes\lambda$ and 
$(\eta_1\otimes\eta_1)\otimes\lambda$
are zero.
Taking $\eta_0$ to be $[1]$ and $\eta_1=[z]$, it follows
that $(\eta_0\otimes\eta_1)\otimes\lambda$ and 
$(\eta_1\otimes\eta_0)\otimes\lambda$ are both 1.
Therefore 
\[
 (\gamma\otimes\gamma)\otimes\lambda=\zeta_0\otimes\zeta_1
  +\zeta_1\otimes\zeta_0.
\]
Taking the spectral triples $(L_2(h),\pi,D)$
and $(\bbc,\sigma,\epsilon,0)$ to represent the classes
$\zeta_1$ and $\zeta_0$ respectively, it follows that the
triple given by $(\clh,\phi,D_0)$ where
\[
 \clh=L_2(h)\oplus L_2(h),\quad
\phi(a\otimes b)=\pi(a)\epsilon(b)\oplus \epsilon(a)\pi(b),\quad
D_0=D\oplus D,
\]
gives the required class. Therefore the restriction of $\phi$ to
the first copy of $\mathcal{A}$ together with $\clh$ and $D$ give a
fundamental class for $\mathcal{A}$.
\qed

\brmrk
{\rm
The $2\ell+1$-dimensional quantum sphere $S_q^{2\ell+1}$
is given by the universal $C^*$-algebra $A_\ell:=C(S_q^{2\ell+1})$ 
generated by
elements
$z_1, z_2,\ldots, z_{\ell+1}$
satisfying the following relations:
\bean
z_i z_j & =& qz_j z_i,\qquad 1\leq j<i\leq \ell+1,\\
z_i^* z_j & =& q z_j z_i^* ,\qquad 1\leq i\neq j\leq \ell+1,\\
z_i z_i^* - z_i^* z_i +
(1-q^{2})\sum_{k>i} z_k z_k^* &=& 0,\qquad \hspace{2em}1\leq i\leq \ell+1,\\
\sum_{i=1}^{\ell+1} z_i z_i^* &=& 1.
\eean
The $K$-theory and the $K$-homology groups for this algebra $A_\ell$
are known and by the same argument as in the proof of theorem~\ref{pd},
$A_\ell$ is $KK$-equivalent to $C(S^1)$ and Poincar\`e duality holds for $A_\ell$.

If one replaces the counit $\epsilon$ for $C(SU_q(2))$ by the functional
\[
z_j\mapsto \cases{1 & if $j=1$,\cr
            0 & if $j\neq 1$,}
\]
on $A_\ell$, and replaces the equivariant spectral triple for $SU_q(2)$
by the spectral triple for $S_q^{2\ell+1}$ equivariant under
the action of $SU_q(\ell+1)$ constructed in~\cite{c-p2}, 
then everything in this section goes through for the odd dimensional quantum spheres.
}
\ermrk


\noindent
{\sc Partha Sarathi Chakraborty}
(\texttt{parthac@imsc.res.in})\\
\begin{footnotesize}School of Mathematical  Sciences, University of Adelaide, AUSTRALIA                                                        \end{footnotesize}\\
         {\footnotesize  (On leave from) Institute of Mathematical Sciences, 
CIT Campus, Chennai--600\,113, INDIA}\\[1ex]
{\sc Arupkumar Pal} (\texttt{arup@isid.ac.in})\\
         {\footnotesize Indian Statistical
Institute, 7, SJSS Marg, New Delhi--110\,016, INDIA}


\end{document}